\documentclass[a4paper,fleqn,11pt]{article}
\usepackage{enumitem}
\usepackage[authoryear]{natbib}
\usepackage[spanish, english]{babel}
\usepackage{moresize}

\usepackage{epic}
\usepackage{calc}
\usepackage{graphicx}
\usepackage{amsmath,amssymb,amsthm,amsfonts}
\usepackage{multirow}
\usepackage{tikz}
\usetikzlibrary{shapes,arrows}
\usetikzlibrary{calc}
\usepackage{fancyhdr, blindtext, textcase}
\usepackage{datetime}
\usepackage{pgffor}
\usepackage{setspace}
\usepackage{pdfpages}
\usepackage{eurosym}
\usepackage{dsfont}
\usepackage{multicol}
\usepackage[normalem]{ulem}
\usepackage{framed}

\usepackage{colortbl}
\usepackage{stackrel}
\usepackage{adjustbox}
\usepackage{empheq}
\usepackage[colorlinks=false, citecolor=blue]{hyperref}

\newcommand{\N}{\mathds{N}}
\newcommand{\R}{\mathds{R}}

\newcommand{\rojo}[1]           {\textcolor[rgb]{0.98,0.00,0.00}{#1}} 
\newcommand{\rojooscuro}[1]     {\textcolor[rgb]{0.90,0.00,0.00}{#1}} 
\newcommand{\azul}[1]           {\textcolor[rgb]{0.00,0.00,1.00}{#1}} 
\newcommand{\naranja}[1]        {\textcolor[RGB]{255,40,0}{#1}} 
\xdefinecolor{steel_blue}{RGB}{35,107,142} 
\xdefinecolor{verdecaja}{RGB}{229,236,228} 
\xdefinecolor{rojocaja}{RGB}{237,218,218} 
\xdefinecolor{azulfondocaja}{RGB}{233,240,244}
\xdefinecolor{azulcaja}{RGB}{222,229,232}
\xdefinecolor{naranja}{RGB}{255,69,0}
\xdefinecolor{azul}{RGB}{0,0,255}
\definecolor{darkolivegreen}{rgb}{0.33, 0.42, 0.18}
\definecolor{olive}{rgb}{0.5, 0.5, 0.0}

\newcommand{\msout}[1]{\text{\sout{\ensuremath{#1}}}}
\newcommand{\rv}[1]{\relax\ifmmode\azul{#1}\else\azul{\textbf{#1}}\fi}
\newcommand{\rvv}[2]{\rojo{\ifmmode\msout{#1}\else\sout{#1}\fi} \azul{\ifmmode#2\else\textbf{#2}\fi}}

\newcommand{\Sebrv}[1]{\relax\ifmmode\naranja{#1}\else\naranja{\textbf{#1}}\fi}
\newcommand{\Sebrvv}[2]{\rojooscuro{\ifmmode\msout{#1}\else\sout{#1}\fi} \naranja{\ifmmode#2\else\textbf{#2}\fi}}

\usepackage[ruled,vlined,linesnumbered]{algorithm2e}
\usepackage[nomessages]{fp}
\usepackage{enumitem}
\usepackage{subcaption}
\usepackage{anysize}
\usepackage{titlesec}
\titlespacing*{\section}{0pt}{0.2cm}{0.1cm}
\titlespacing*{\subsection}{0pt}{0.2cm}{0.1cm}
\titlespacing*{\title}{0pt}{0.3cm}{0.3cm}
\titlespacing*{\abstract}{0pt}{0.1cm}{0.3cm}
\marginsize{2.1cm}{2.1cm}{0.9cm}{0.9cm}
\makeatletter
\newcommand\subsubsubsection{\@startsection{paragraph}{4}{\z@}{-2.5ex\@plus -1ex \@minus -.25ex}{1.25ex \@plus .25ex}{\normalfont\normalsize\bfseries}}
\newcommand\subsubsubsubsection{\@startsection{subparagraph}{5}{\z@}{-2.5ex\@plus -1ex \@minus -.25ex}{1.25ex \@plus .25ex}{\normalfont\normalsize\bfseries}}
\makeatother

\usepackage[utf8]{inputenc}
\usepackage[T1]{fontenc}
\usepackage[spanish]{babel}

\usepackage{array}
\newcolumntype{C}[1]{>{\centering\arraybackslash}m{#1}}

\title{A Survey of Community Detection from an Operations Research Perspective: Taxonomy, Mathematical Formulations, Modularity Functions, and Benchmark Datasets}
\author{
  Pozo Montaño, Miguel A.        \footnote{Department of Statistics and Operational Research \& Institute of Mathematics (IMUS), University of Seville. E-mail: miguelpozo@us.es.}\\
  Taboh, Sebastián V.       \footnote{Corresponding author. Department of Computer Science, University of Buenos Aires - Institute of Computer Science (UBA/CONICET). E-mail: staboh@dc.uba.ar.}
}
\date{\today}

\newcommand{\Cross}{\includegraphics[width = 0.3cm, height = 0.25cm]{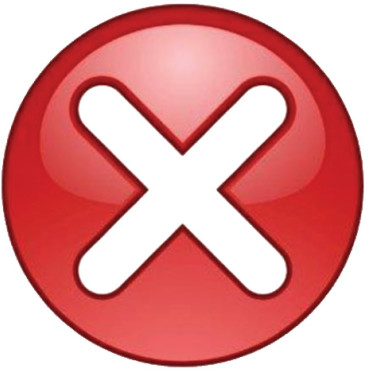}}

\newcommand{\Checkmark}{\centering \includegraphics[width = 0.3cm, height = 0.25cm]{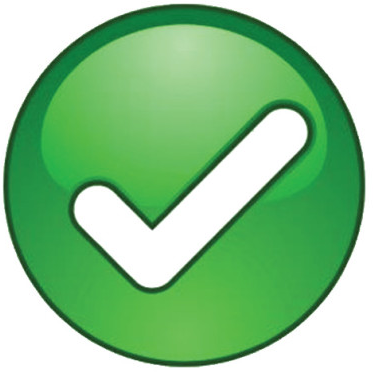}}

\usepackage{makecell}

\newtheorem{definition}{Definition}[section]

\newcommand{\Soc}{S}

\newcommand{\BioEco}{BEHS}
\newcommand{\Comm}{CT}
\newcommand{\EComm}{EE}
\newcommand{\Aca}{AS}
\newcommand{\Info}{IS}

\usepackage{booktabs}

\newtoggle{tablaConReferencias}
\togglefalse{tablaConReferencias}
\newcommand{\tabcite}[1]{%
  \iftoggle{tablaConReferencias}{~\cite{#1}}{}%
}

\begin{document}

\maketitle

\begin{abstract}
Community detection is a fundamental problem in network science that consists of identifying groups of vertices exhibiting stronger internal connectivity than external connectivity. From an Operations Research perspective, the problem can be interpreted as a family of combinatorial optimization and clustering models defined over network structures. This survey provides a unified framework for the study of community detection, with particular emphasis on modularity-based approaches. We first review existing surveys and analyze the taxonomic criteria used to classify the literature, highlighting the absence of a common conceptual framework. Based on this analysis, we propose a multidimensional taxonomy that organizes community detection methods according to network characteristics, community structure, objective functions, methodological paradigms, evaluation criteria, and application domains. We then introduce a general mathematical formalization of the Community Detection Problem that accommodates disjoint, overlapping, and fuzzy community structures within a unified assignment framework. Building on this formalization, we review representative modularity functions, discussing their underlying assumptions, null models, and known limitations. We also survey modularity-based community detection methods, distinguishing between algorithmic and mathematical programming approaches. Finally, we review commonly used benchmark datasets and discuss their role in evaluation and reproducibility. By integrating taxonomy, mathematical modeling, modularity analysis, and benchmarking practices, this survey provides a structured reference for researchers and practitioners working on community detection and related network optimization problems.


\end{abstract}

\section{Introduction}
\label{sec:Introduction}

Graphs constitute one of the most powerful abstract structures for modeling relationships among interacting entities. They naturally arise in a wide variety of contexts, including social systems, communication networks, biological interactions, financial systems, and technological infrastructures. In such settings, vertices represent entities and edges encode interactions or relationships between them. As the scale and complexity of relational data have increased, so has the need for systematic methods capable of uncovering latent structural patterns within networks.

One of the most prominent structural patterns studied in network analysis is the presence of communities, loosely defined as groups of vertices that are more densely connected internally than externally. The \emph{Community Detection Problem} consists of identifying such groups in a given graph. Although conceptually intuitive, this problem admits multiple formal interpretations and modeling choices. Depending on the application context, networks may be directed or undirected, weighted or unweighted, static or evolving. Likewise, communities may be required to form a partition of the vertex set (disjoint communities) or may overlap, allowing vertices to belong to multiple groups. These modeling alternatives significantly affect both the mathematical formulation of the problem and the methodological approaches used to solve it.

From the perspective of Operations Research and Management Science (OR/MS), the Community Detection Problem can be naturally interpreted as a family of optimization and clustering problems. At its core, the task consists of assigning vertices to one or more subsets so as to optimize a structural objective function that captures the quality of the resulting partition or cover. Among the various objective functions proposed in the literature, modularity has played a central role since its introduction \citep{Newman2004}, becoming one of the most widely adopted criteria for evaluating community structure. However, despite its popularity, modularity and its variants have been defined and employed under heterogeneous assumptions, often without a unified notation or explicit discussion of their modeling implications.

The rapid growth of the community detection literature has led to a highly fragmented landscape. Numerous surveys have been published across disciplines such as physics, computer science, applied mathematics, and social network analysis. These reviews differ substantially in terminology, taxonomic criteria, problem formalization, and evaluation perspectives. In many cases, the problem is described procedurally (focusing on algorithms) rather than structurally (emphasizing the underlying optimization models and assumptions). As a result, there is limited consensus regarding how to systematically classify methods, how to formally compare alternative formulations of modularity, or how to relate algorithmic approaches to precise mathematical models. This fragmentation poses particular challenges for the OR/MS community. Without a coherent conceptual framework, it becomes difficult to interpret community detection as a well-defined class of optimization problems, to understand the relationships among existing models, or to identify meaningful research gaps. Moreover, inconsistencies in definitions and evaluation criteria block reproducibility and complicate the assessment of methodological advances.

The primary objective of this paper is to provide a unifying conceptual and methodological framework for the Community Detection Problem from an OR/MS perspective. Rather than proposing new algorithms or computational experiments, we aim to clarify, systematize, and formalize the existing body of knowledge. The contribution of this review is therefore to organize a mature but conceptually dispersed field into a coherent framework that facilitates modeling, comparison, and future development.

The paper makes four main contributions.

First, we conduct a comprehensive review of existing surveys in community detection. This meta-review highlights the absence of consensus in taxonomy, terminology, and problem framing, thereby motivating the need for a new classification scheme explicitly designed to support modeling and decision-oriented analysis.

Second, we propose a taxonomy of tags that organizes the literature along structural, methodological, and operational dimensions. The taxonomy distinguishes, among other aspects, network characteristics (e.g., directed vs.\ undirected, weighted vs.\ unweighted), community structures (disjoint vs.\ overlapping), methodological paradigms (exact, heuristic, probabilistic, and learning-based approaches), and practical challenges (e.g., scalability, resolution limits, stability). This taxonomy is not merely descriptive; it is intended as a conceptual tool to support the formulation and comparison of models within an OR/MS framework.

Third, we provide a unified formalization of the Community Detection Problem and of modularity-based objective functions. By consolidating notation, modeling assumptions (e.g., disjoint versus overlapping communities, fuzzy memberships, and network variants), and modularity formulations scattered across the literature, we clarify how commonly used formulations can be interpreted within a common assignment-based framework and identify the implicit modeling choices underlying widely used methods. Although no new theoretical results are introduced, this normalization reduces ambiguity and enhances the interpretability and reproducibility of future research.

Fourth, we review the principal datasets employed in community detection studies, interpreting them not simply as empirical examples but as benchmark instances that influence how algorithms are evaluated, compared, and validated. This perspective is consistent with OR/MS traditions, where data play a central role in validating and comparing methodological approaches.

By positioning community detection explicitly as a family of optimization and clustering problems, this review seeks to bridge the gap between the interdisciplinary network science literature and the OR/MS community. The framework developed here provides a common language and a structured view of the field, supporting both theoretical modeling and applied research.

The remainder of the paper is organized as follows. 
Section 2 presents a comprehensive review of existing surveys and discusses their scope, taxonomic criteria, and methodological emphasis. 
Section 3 introduces the proposed taxonomy and explains its structural, methodological, and operational dimensions.
Section 4 formalizes consistent definitions of communities, assignments of elements into communities, modularity, and the community detection problem.
Section 5 reviews modularity-based community detection methods, distinguishing between exact and heuristic approaches, and highlighting formulations grounded in mathematical programming. 
Section 6 presents and discusses representative modularity functions used in the literature.
Section 7 reviews the datasets most frequently used in the literature and discusses their role in benchmarking and evaluation. 
Finally, Section 8 concludes the paper and outlines open issues and directions for future research.

\section{Literature Review on the Community Detection Problem}

Community detection is one of the central problems in network science, concerned with identifying groups of nodes that are more densely connected internally than with the rest of the network. Over the past two decades, this problem has attracted considerable attention due to its importance in a wide range of domains, including social networks, biological systems, information networks, and infrastructure analysis. As a result, an extensive literature has emerged, proposing diverse models and algorithms for detecting meaningful communities in complex networks.

Among the various methodological families, a prominent group of approaches is modularity-based, in which the goal is to optimize a modularity function that quantifies the quality of a partition by comparing the observed density of edges within communities against that expected under a suitable null model. High modularity values indicate a significant deviation from randomness, suggesting well-defined community structures. In contrast, many other approaches do not rely on modularity optimization but instead stem from alternative principles, such as spectral clustering, probabilistic modeling, information-theoretic formulations, or learning-based paradigms.

To synthesize this wide range of developments, several surveys have been published that classify and compare community detection methods from different perspectives. Some organize the literature according to the underlying graph or learning models, while others focus on algorithmic strategies, optimization objectives, or types of community structures. Together, these surveys provide a comprehensive view of the evolution and diversity of approaches to community detection.

The following paragraphs summarize the main reviews on the Community Detection Problem published between 2013 and 2024.
The emphasis is placed on the scope of each review, the taxonomic dimensions (or classification criteria) adopted to classify the literature, and any distinctive contributions such as benchmark analyses, datasets, or emerging research directions.
The reviews are discussed in chronological order to highlight the evolution of taxonomic perspectives over time.

\subsection{Early Reviews Focused on Structural Aspects (2013--2015)}

\citet{Xie2013} represents one of the earliest comprehensive surveys dedicated exclusively to the problem of overlapping community detection.
Rather than addressing community detection in general, the review deliberately narrows its scope to networks where nodes may belong to multiple communities simultaneously.
The classification of methods is primarily algorithmic, grouping approaches according to their underlying detection strategy.
Major categories include clique percolation methods, seed expansion techniques, link-based approaches, and probabilistic or fuzzy models.
The taxonomy is essentially single-axis, focused on methodological principles rather than network structure or application domain.
A distinctive feature of this review is its strong empirical orientation.
The authors provide comparative experiments using synthetic benchmark networks and discuss performance differences across methods.
This work establishes overlapping membership as a fundamental taxonomic dimension, later reused by many subsequent surveys.

\citet{SurveyBouhali2015} focuses specifically on community detection in social networks, with an emphasis on classical methods and open research challenges.
The scope is narrower than later surveys, concentrating on social interaction graphs rather than general complex networks.
The classification is largely method-oriented, distinguishing between hierarchical clustering approaches, partition-based methods, spectral techniques, and modularity-based algorithms.
Structural properties such as directed or weighted networks are occasionally mentioned but do not form explicit taxonomic axes.
The review is mostly descriptive and does not propose a formal or multi-dimensional taxonomy.
Its main contribution lies in identifying research perspectives and unresolved issues, such as scalability and dynamic behavior, which are discussed qualitatively rather than systematically evaluated.

\subsection{Conceptual Consolidation and Multidimensional Views (2016--2018)}

The highly influential review of \citet{Fortunato2016} adopts a conceptual and explanatory perspective rather than aiming for exhaustive algorithmic coverage.
Its primary objective is to clarify what is meant by a ``community'' and to explain why different algorithms often yield different partitions for the same network.
Instead of proposing a strict taxonomy, the authors organize the discussion around several conceptual dimensions:
definitions of communities (e.g., density-based, flow-based, similarity-based),
objective functions (modularity, stability, likelihood),
and known limitations such as the resolution limit and solution degeneracy.
These dimensions implicitly act as classification axes.
A key contribution of this work is the explicit integration of challenges and limitations into the discussion.
Rather than treating them as secondary issues, they are presented as intrinsic to the problem itself.
This perspective strongly motivates later attempts to formalize challenges as a separate taxonomic category.

\citet{Schaub2017} further develops the idea that community detection is inherently multi-faceted and context-dependent.
The authors argue that no single definition or algorithm can be universally optimal.
The organization of the review emphasizes multiple perspectives, including the role of the objective function, the scale or resolution of communities, and the distinction between static and dynamic networks.
Rather than proposing a hierarchical taxonomy, the paper promotes a multidimensional conceptual framework.
An important contribution is the explicit recognition that different taxonomies may be equally valid depending on the research question.
This view directly supports the need for flexible, multi-axis classification schemes, such as the taxonomy proposed in the present work.

\citet{Javed2018} aims to provide a broad overview of community detection across multiple disciplines, including computer science, physics, and social sciences.
The taxonomy is explicitly hierarchical.
At the top level, methods are classified according to whether they detect disjoint or overlapping communities.
Within each category, further classification is performed based on the algorithmic approach, such as modularity-based, label propagation, or dynamic methods.
Dynamic networks are treated as a separate category.
This review represents one of the first systematic attempts to structure the literature using clearly defined hierarchical tags (that is, categories or classification levels).
It highlights the importance of the type of community as a primary taxonomic axis.

\subsection{Expansion Toward Applications and Specialized Network Types (2019--2021)}

\citet{SurveyElMoussaoui2019} provide a broad overview of both methodological approaches and application domains.
Their scope includes social, biological, and technological networks.
The classification combines algorithmic criteria with application-oriented groupings.
Methods are categorized according to their general approach (e.g., modularity-based, spectral, heuristic), while applications are discussed separately.
The taxonomy is therefore hybrid rather than strictly hierarchical.
A notable feature is the emphasis on practical applications, which reinforces the relevance of the application domain as an important taxonomic dimension.

\citet{Huang2021} focus exclusively on multilayer and multiplex networks, which involve multiple types of interactions or relational contexts.
The taxonomy is centered on the structure of the network itself, distinguishing between different multilayer models and corresponding detection methods.
Methods are classified according to how they handle inter-layer dependencies and cross-layer consistency.
This work clearly justifies treating multilayer structure as a primary classification axis rather than a minor extension of single-layer networks.

\citet{SurveyJin2021} capture the methodological transition from classical statistical models to deep learning-based approaches.
The classification is organized by methodological paradigm, distinguishing between probabilistic generative models, traditional graph-based algorithms, and deep learning methods such as autoencoders and graph neural networks.
The review highlights the growing importance of representation learning and motivates the inclusion of the learning paradigm as a key taxonomic dimension.

\citet{Attea2021} present a comprehensive survey of heuristic- and metaheuristic-based community detection algorithms. The authors also introduce two new taxonomies  --- hybrid metaheuristic (HMCD) and hyper heuristic (HHCD) frameworks --- as design paradigms for future community detection algorithms that more effectively exploit problem-specific knowledge through cross-fertilization between heuristic and metaheuristic components.

\subsection{Recent Specialization and Fragmentation (2022--2024)}

\citet{SurveyAlotaibi2022} concentrate on temporal and dynamic networks, particularly in social contexts.
The primary classification axis is temporality, distinguishing between static snapshots and explicitly time-evolving models.
Within each category, methods are further classified by algorithmic strategy.
This review reinforces temporality as a fundamental structural dimension rather than a secondary extension.

\citet{Su2024} provide an in-depth survey of deep learning approaches to community detection.
The taxonomy is highly structured, with an initial division between traditional methods and deep learning-based approaches.
Deep learning methods are further subdivided into autoencoder-based models, graph neural networks, reinforcement learning approaches, and hybrid architectures.
While the taxonomy is detailed, its scope is restricted to a specific methodological paradigm.
This review clearly isolates the learning paradigm as an independent classification axis, supporting its inclusion as a primary dimension in modern taxonomies of community detection methods.

\citet{Diboune2024} propose a distinctive classification based on the nature of the contribution.
At the top level, works are divided into theoretical and applied contributions.
Theoretical approaches are further classified according to the underlying model, including graph-based, machine learning, multi-objective optimization, and game-theoretic models.
This contribution-oriented taxonomy departs from purely algorithmic classifications and introduces a novel perspective.
This perspective suggests that taxonomies can also be defined at the level of research contribution, rather than purely methodological or structural criteria, 
introducing a complementary dimension that is rarely formalized in earlier surveys.

\citet{Li2024} offer a modern, general overview of community detection methods.
The classification blends multiple criteria, including network structure, algorithmic approach, and evaluation methodology.
While comprehensive, the taxonomy is less hierarchical and serves more as a mapping of the research landscape.
The absence of a strict hierarchical structure in this taxonomy further illustrates the lack of consensus in the literature, 
and motivates the need for a more systematic and unified framework.

\citet{Tekin2024} address adversarial attacks on community detection algorithms.
The taxonomy is centered on the type of attack and threat model rather than on detection methods themselves.
Although orthogonal to classical taxonomies, this review introduces security and robustness as transversal concerns.

\subsection{Synthesis and Motivation for a Unified Taxonomy}

The reviewed surveys demonstrate a wide diversity of taxonomic criteria, including algorithmic strategy, network structure, type of community, learning paradigm, application domain, and contribution type.
No single taxonomy dominates the literature, and most reviews adopt a partial or context-specific perspective.
This fragmentation motivates the need for a unified, multi-dimensional taxonomy capable of integrating these heterogeneous viewpoints.

\section{Taxonomy of Classification Tags for the Community Detection Problem}
	
The field of community detection in complex networks has grown remarkably diverse, encompassing a vast range of algorithms, objectives, and application contexts. However, the absence of a unified conceptual framework often leads to fragmented methodological comparisons and terminological inconsistencies across studies. To address this issue, we propose a structured taxonomy that organizes existing research into six complementary dimensions: 
(1) \textit{Network Structure}, describing the structural characteristics of the analyzed graphs; 
(2) \textit{Community Structure}, defining the nature and membership of the groups sought; 
(3) \textit{Objective or Quality Function}, specifying the criteria optimized by detection methods; 
(4) \textit{Evaluation and Validation}, covering how community structures are assessed and benchmarked;
(5) \textit{Application Domain}, identifying the practical contexts where community detection is applied and
(6) \textit{Method or Solution Technique}, detailing the algorithmic strategies employed.
This taxonomy provides a coherent framework to compare methods, reveal conceptual relationships, and to guide future research toward addressing persistent challenges related to scalability, interpretability, and robustness.

The taxonomy proposed in this section is directly informed by the analysis of existing surveys presented in Section 2. 
Specifically, each dimension of the taxonomy reflects recurring classification criteria identified across the literature. 
For instance, structural properties such as directionality, multilayer structure, and temporality emerge from specialized surveys on network types (e.g., \citealp{Huang2021}; \citealp{SurveyAlotaibi2022}), 
while distinctions between disjoint and overlapping communities are consistently emphasized in early and multidisciplinary reviews (e.g., \citealp{Xie2013}; \citealp{Javed2018}). 
Similarly, methodological categorizations (from classical optimization techniques to deep learning approaches) reflect the evolution highlighted in more recent surveys (e.g., \citealp{SurveyJin2021}; \citealp{Su2024}). 

Overall, the proposed taxonomy can be interpreted as a synthesis and formalization of these heterogeneous perspectives, reorganizing them into a coherent multi-dimensional framework. The following subsections describe each dimension and its associated classification tags in detail.

\subsection{Network structure}
These tags describe the structural properties of the type of network being analyzed.

\begin{itemize}[leftmargin=0.5cm]
	
	\item \textbf{Directionality}: (Directed $|$ Undirected)  
	Indicates whether the network edges have inherent direction.  
	\emph{Directed} networks represent asymmetric relationships (e.g., follower--followee links on Twitter),  
	while \emph{Undirected} networks assume mutual connections (e.g., friendships).  
	Directionality strongly affects the definition and interpretation of communities, as direction may encode hierarchy or influence.
	
	\item \textbf{Weight}: (Weighted $|$ Unweighted) 
	Specifies whether edges carry numerical values representing intensity, frequency, or capacity.  
	\emph{Weighted} networks encode additional relational information (e.g., number of interactions),  
	whereas \emph{Unweighted} ones consider only presence or absence of links.  
	Some algorithms can generalize modularity or similarity measures to incorporate weights.
	
	\item \textbf{Layers}: (Simple $|$ Multilayer)  
	Refers to the number of distinct relationship types represented in the network.  
	A \emph{Simple} network includes a single connection type;  
	\emph{Multilayer} networks contain several layers,  
	where each layer corresponds to a different interaction context (e.g., social, economic, transportation).  
	These structures allow community detection across layers or within individual ones.
	
	\item \textbf{Attributed networks}: (Node-attributed $|$ Edge-attributed)  
	Indicates the presence of metadata associated with nodes or edges.  
	Node attributes can include demographic, spatial, or textual information,  
	while edge attributes may represent interaction type or context.  
	Community detection in attributed networks combines topological and semantic similarity.
	
	\item \textbf{Hypergraph structure}: (Simple graph $|$ Hypergraph)  
	Extends the graph model so that edges (hyperedges) can connect more than two nodes simultaneously.  
	\emph{Hypergraphs} capture group interactions (e.g., co-authorship or collaboration networks).  
	Community detection methods for hypergraphs typically rely on tensor decompositions or hyperedge clustering.
	
	\item \textbf{Temporality}: (Static $|$ Dynamic)  
	Describes whether the network structure evolves over time.  
	\emph{Static} networks are fixed and analyzed as a single snapshot.  
	\emph{Dynamic} or \emph{Temporal} networks explicitly consider time-evolving links or communities,  
	often using streaming or time-slice models.  
	Temporal dynamics may reveal community formation, persistence, and dissolution patterns.
	
	\item \textbf{Heterogeneity}: (Homogeneous $|$ Heterogeneous $|$ Multitype)  
	Differentiates networks with identical node/edge types (\emph{Homogeneous})  
	from those containing multiple types or relations (\emph{Heterogeneous}).  
	Heterogeneous or multitype networks include, for example, user–item or author–paper graphs.  
	Detecting communities in such networks may involve meta-path or co-clustering approaches.
	
	\item \textbf{Partition structure}: (Bipartite $|$ Multipartite)  
	Refers to the inherent topological division of the input nodes into disjoint sets, where edges occur only between (not within) sets.  
	\emph{Bipartite} graphs connect two sets (e.g., user–movie networks),  
	while \emph{Multipartite} generalizes to more than two disjoint sets.
	Community detection here often requires projection or tailored modularity definitions.
	
\end{itemize}

\subsubsection{Note on conceptual overlap.}  
The categories \emph{Layers}, \emph{Attributed networks}, and \emph{Hypergraph structure} partially overlap,  
as they all extend the classical single-layer, homogeneous graph model by introducing richer relational contexts.  
In general, \emph{multilayer networks} represent multiple relation types (e.g., social vs.\ professional),  
\emph{attributed networks} enrich nodes or edges with auxiliary information (e.g., demographics or text),  
and \emph{hypergraphs} capture higher-order interactions among groups of nodes.  
\emph{Heterogeneity} further generalizes these ideas by allowing multiple node or edge types,  
and can conceptually subsume some multilayer or attributed formulations depending on representation.  
Despite these intersections, they are treated as distinct dimensions for analytical clarity, as each has led to specific  
community detection formulations and methodological frameworks in the literature.

\subsection{Community structure}
These tags indicate the structural nature and membership properties of the detected communities or clusters.

\begin{itemize}[leftmargin=0.5cm]
	
	\item \textbf{Intersections}: (Disjoint $|$ Overlapping $|$ Hybrid)  
	Defines whether nodes can belong to more than one community.  
	\emph{Disjoint} (or non-overlapping) communities assign each node to exactly one group, producing a clear partition.  
	\emph{Overlapping} communities allow nodes to participate in multiple groups, reflecting realistic social or biological settings.  
	\emph{Hybrid} approaches can detect both types depending on network characteristics or parameters.
	
	\item \textbf{Hierarchicity}: (Flat $|$ Hierarchical $|$ Multilevel)  
	Indicates whether communities are organized across multiple levels of granularity.  
	\emph{Flat} detection yields a single partition,  
	while \emph{Hierarchical} (or \emph{Nested}) methods uncover sub-communities within larger ones, often visualized as dendrograms.  
	\emph{Multilevel} approaches allow zooming between coarse and fine resolutions, integrating different scales.
	
	\item \textbf{Membership}: (Hard $|$ Fuzzy $|$ Probabilistic)  
	Specifies how membership is defined.  
	\emph{Hard} clustering assumes binary membership (belong / not belong).  
	\emph{Fuzzy} or \emph{Soft} clustering assigns degrees of membership to each node,  
	while \emph{Probabilistic} methods interpret these degrees as probabilities under a generative model.  
	Fuzzy memberships are common in latent space or NMF-based community models.
	
	\item \textbf{Community unit}: (Node-based $|$ Edge-based)  
	In \emph{Node-based} approaches, communities group vertices with dense interconnections.
    \emph{Edge-based} approaches cluster edges instead, producing overlapping node communities implicitly.
	Edge-centric methods are especially useful for analyzing interaction patterns or role-based networks.
	
	\item \textbf{Scale}: (Single-scale $|$ Multi-scale $|$ Resolution-adaptive)  
	Refers to whether communities are detected at a single resolution or across multiple levels of granularity.  
	\emph{Single-scale} methods use a fixed resolution parameter.  
	\emph{Multi-scale} methods detect communities of varying sizes simultaneously,  
	and \emph{Resolution-adaptive} methods automatically adjust the detection scale to network structure.
	
\end{itemize}

\subsection{Objective or Quality Function}
These tags describe the criteria used to evaluate or optimize the community partition quality.  
They define the mathematical or conceptual basis guiding the evaluation or the optimization process.

\begin{itemize}[leftmargin=0.5cm]
	
	\item \textbf{Modularity-based}: (Standard modularity $|$ Weighted modularity $|$ Multi-resolution modularity $|$ Signed modularity)  
	The most widely used measure; it compares observed intra-community connections to those expected under a null (random) model.  
	Extensions incorporate weights, resolution parameters, or signed edges.

	\item \textbf{Likelihood / Statistical-model based}: (Stochastic Block Model $|$ Degree-corrected SBM $|$ Mixed-membership SBM $|$ Bayesian inference)  
	Assumes that the observed network is generated by a probabilistic model where communities correspond to latent variables.  
	Community detection then becomes a maximum-likelihood or Bayesian inference problem.  
	These approaches provide interpretability and statistical guarantees.
	
	\item \textbf{Random-walk / Flow-based}: (Infomap $|$ Markov stability $|$ Flow retention models)  
	Defines communities based on the tendency of random walks to remain within the same group for long periods.  
	\emph{Infomap} uses information-theoretic compression of random walk trajectories.  
	\textit{Flow-based} approaches capture functional modularity, relevant for transport or communication systems.
	
	\item \textbf{Information-theoretic}: (Entropy minimization $|$ Description length $|$ Mutual information)  
	Measures community quality by minimizing information loss or description length when representing the network structure.  
	These approaches interpret communities as efficient encodings or compressions of connectivity data.
	
	\item \textbf{Cut / Partition metrics}: (Min-cut $|$ Ratio-cut $|$ Normalized-cut $|$ Conductance $|$ Edge-betweenness)  
	Aim to minimize the number or weight of inter-community edges relative to intra-community ones.  
	\emph{Min-cut} and its normalized variants are classical in spectral clustering and graph partitioning.  
	\textit{Edge-betweenness} methods remove highly ``bridging'' edges iteratively to reveal modular structures.
	
	\item \textbf{Stability / Persistence}: (Temporal persistence $|$ Perturbation stability $|$ Ensemble stability)  
	Assess whether communities remain consistent over time, under noise, or across different runs of the algorithm.  
	\emph{Temporal persistence} applies to evolving networks,  
	while \emph{Perturbation stability} evaluates robustness to small topological changes.  
	\textit{Ensemble stability} measures consistency across multiple algorithmic realizations.
		
	\item \textbf{Multi-objective}: (Weighted-sum $|$ Pareto-optimal $|$ Hybrid)  
	Combines several quality functions (e.g., modularity, conductance, attribute similarity) into a single optimization problem.  
	\emph{Weighted-sum} methods aggregate objectives linearly,
    while \emph{Pareto-optimal} approaches search for non-dominated trade-offs among competing criteria.
    Hybrid approaches combine aggregation and Pareto-based mechanisms.
	
\end{itemize}

\subsection{Evaluation and Validation}
These tags refer to how the quality and reliability of community detection results are assessed.  
They determine whether the detected structure is meaningful, robust, and comparable across algorithms or datasets.

\begin{itemize}[leftmargin=0.5cm]
	
	\item \textbf{Ground truth / Benchmark datasets}: (Available $|$ Unavailable $|$ Partial labels)  
	Indicates whether known community labels exist for validation.  
	\emph{Available} ground truth allows direct comparison,
	while \emph{Unavailable} cases require unsupervised or intrinsic metrics.  
	\emph{Partial labels} indicate that community information is available only for a subset of vertices or communities.
	
	\item \textbf{Data realism}: (Synthetic $|$ Real-world $|$ Hybrid)  
	\emph{Synthetic} networks are generated with known ground-truth communities for controlled testing and scalability evaluation.  
	\emph{Real-world} networks assess applicability and interpretability under natural noise and heterogeneity.  
	\emph{Hybrid} studies combine both to validate generalization and realism.
	
	\item \textbf{Evaluation metrics}: (External $|$ Internal $|$ Relative $|$ Computational)  
	\emph{External} metrics compare detected communities with known labels.  
	\emph{Internal} metrics assess structure quality without ground truth.  
	\emph{Relative} metrics benchmark performance against baselines or null models. 
	\emph{Computational} metrics include runtime, scalability, and memory efficiency.
	
	\item \textbf{Consistency (stability under perturbations)}: (Node/edge removal $|$ Edge rewiring $|$  Noise injection)  
	Evaluates how consistent the detected communities remain under small structural changes.
	Methods may test robustness against \emph{node/edge removal}, random \emph{rewiring}, or controlled \emph{noise injection}.
	
	\item \textbf{Robustness (attack tolerance)}: (Noise robustness $|$ Adversarial tolerance $|$ Missing data resilience)  
	Measures algorithmic sensitivity to incomplete, noisy, or manipulated data.  
	\emph{Noise robustness} evaluates random errors;  
	\emph{Adversarial tolerance} studies resilience against targeted attacks or misinformation;  
	\emph{Missing data resilience} measures the ability to recover structure when edges or nodes are absent.
	
	\item \textbf{Scalability/Architecture}: (Small-scale $|$ Large-scale $|$ Parallel $|$ Distributed)  
	Examines computational efficiency as network size grows.  
	\emph{Small-scale} algorithms suit networks up to a few thousand nodes.  
	\emph{Large-scale} methods exploit heuristics or approximations to maintain viability on larger networks.  
	\emph{Parallelizable} and \emph{Distributed} approaches leverage multi-core or cluster architectures respectively, to scale further.
	
	\item \textbf{Interpretability}: (Transparent $|$ Black-box $|$ Explainable)  
	Refers to how easily the detected communities and algorithmic reasoning can be understood by humans.  
	\emph{Transparent} models allow direct interpretation of structure.  
	\emph{Black-box} methods require post-hoc \emph{Explainable AI} tools to extract insights.
	
	\item \textbf{Statistical significance}: ($p$-value $|$ Null model comparison $|$ Randomization test)  
	Quantifies whether the detected communities are unlikely to appear by random chance.  
	Common approaches use \emph{p-values}, \emph{Z-scores}, or comparisons with randomized networks.  
	This validation ensures that observed modular structures reflect genuine organization rather than statistical artifacts.
	
\end{itemize}

\subsection{Application Domain}
These tags describe the type of real-world systems where community detection methods are applied.  
Each domain provides unique structural patterns, data types, and interpretative goals.

\begin{itemize}[leftmargin=0.5cm]
	
	\item \textbf{Social networks}: (Friendship $|$ Communication $|$ Collaboration $|$ Online platforms)  
	Networks where nodes represent individuals or groups and edges encode social interactions.  
	Examples include \textit{friendship} graphs, email exchanges, or co-authorship networks.  
	Communities often correspond to social circles, teams, or interest groups.
	
	\item \textbf{Biological networks}: (Protein interaction $|$ Gene co-expression $|$ Metabolic $|$ Ecological)  
	Represent molecular or ecological relationships among biological entities.  
	Communities may represent \textit{protein} complexes, \textit{gene} modules, \textit{metabolic} pathways, or \textit{ecological} niches.  
	These analyses help infer functional organization and biological processes.
	
	\item \textbf{Technological networks}: (Internet $|$ Power grid $|$ Communication infrastructure)  
	Capture connections among technological components such as routers, power stations, or devices.  
	Community detection here aids fault localization, resilience analysis, and functional decomposition.
	
	\item \textbf{Information networks}: (Citation $|$ Co-authorship $|$ Web hyperlink $|$ Knowledge graphs)  
	Represent information flow or association patterns between documents, authors, or webpages.  
	Communities identify topical clusters, scientific disciplines, or linked knowledge domains.
	
	\item \textbf{Transportation networks}: (Air traffic $|$ Road $|$ Maritime $|$ Urban mobility)  
	Nodes represent stations, airports, or intersections; edges denote routes or travel paths.  
	Community detection uncovers modular transport zones or travel corridors, supporting urban and logistics planning.
	
	\item \textbf{Financial / Economic networks}: (Trade $|$ Interbank $|$ Stock correlation $|$ Transaction)  
	Capture dependencies or exchanges among economic entities.  
	Communities correspond to market sectors, trade blocs, or clusters of correlated financial instruments,  
	providing insight into systemic risk and market structure.
	
	\item \textbf{Attribute-rich networks}: (Textual $|$ Categorical $|$ Spatial $|$ Multimodal)  
	Networks where nodes or edges are enriched with descriptive attributes (e.g., user profiles, document topics).  
	Detection integrates both topological proximity and attribute similarity,  
	often through embedding or joint optimization frameworks.
	
	\item \textbf{Spatio-temporal networks}: (Spatial $|$ Temporal $|$ Mobility)  
	Combine spatial and temporal dimensions, such as transportation or epidemic networks.  
	Communities may represent geographically or temporally cohesive groups,  
	revealing migration flows, diffusion patterns, or synchronized activities.
	
\end{itemize}

\subsection{Method or Solution Technique}
This tag categorizes community detection methods according to the underlying computational or theoretical approach used to identify communities. Each main family groups techniques sharing similar computational principles, modeling assumptions or optimization goals.

\begin{itemize}[leftmargin=0.5cm]
	
	\item \textbf{Traditional (or classical) methods}: (Spectral methods $|$ Hierarchical clustering $|$ Graph partitioning $|$ Label propagation)
		\emph{Spectral methods} use eigen-decomposition of adjacency or Laplacian matrices to cluster nodes based on low-dimensional representations.
		\emph{Hierarchical clustering} builds a dendrogram of nested communities via agglomerative or divisive strategies.
		\emph{Graph partitioning} minimizes edge cuts between subsets, e.g., ratio-cut or normalized-cut criteria.
		\emph{Label propagation} iteratively updates node labels based on neighborhood majority, converging to stable community assignments.

	\item \textbf{Optimization-based methods}: (Modularity maximization $|$ Mathematical programming $|$ Metaheuristics)
		\emph{Modularity maximization} searches for partitions that maximize modularity or similar quality functions.
		\emph{Mathematical programming} formulates the problem as MILP, SDP, or other combinatorial optimization models solved exactly or approximately.      \emph{Metaheuristics} leverage global optimization strategies (such as genetic algorithms, simulated annealing, tabu search, or swarm intelligence) to optimize modularity or related measures.

	\item \textbf{Statistical / Generative models}: (Stochastic block models $|$ Degree-corrected SBMs $|$ Mixed-membership models $|$ Bayesian inference)
		\emph{Stochastic Block Models (SBM)} assume edges are generated probabilistically depending on community memberships.
		\emph{Degree-corrected SBM} incorporates node degree heterogeneity to better model real networks.
		\emph{Mixed-membership models} allow nodes to belong to multiple communities with probabilistic weights.
		\emph{Bayesian inference} estimates posterior distributions over partitions given observed data and priors.

	\item \textbf{Graph representation and embedding methods}: (Matrix factorization $|$ Graph embeddings $|$ Link clustering)
		\emph{Matrix factorization} uses NMF, SVD, or tensor decomposition to reveal latent community structures.
		\emph{Graph embeddings} learn low-dimensional vector representations of nodes preserving structural proximity, followed by clustering (e.g., $k$-means).
		\emph{Link clustering} detects overlapping communities by clustering edges rather than nodes.
	
	\item \textbf{Deep learning methods}: (Graph Neural networks $|$ Autoencoders / Variational autoencoders $|$ Generative adversarial networks)
		\emph{Graph Neural Networks (GNNs)} learn hierarchical node representations using message passing schemes, often optimized for clustering loss functions.
		\emph{Autoencoders / VAEs} encode adjacency information into latent embeddings and reconstruct to infer community structure.
		\emph{Generative Adversarial Networks (GANs)} model the distribution of node embeddings or graphs to generate and detect communities.

	\item \textbf{Hybrid and ensemble methods}: (Consensus clustering $|$ Multi-objective optimization $|$ Attribute-structure integration)
		\emph{Consensus clustering} combines results from multiple algorithms to obtain a stable consensus partition.
		\emph{Multi-objective optimization} balances several quality functions (e.g., modularity and density) simultaneously.
		\emph{Attribute-structure integration} jointly exploits topological and attribute information to improve detection accuracy.

	\item \textbf{Dynamic, Streaming, and Multi-layer methods}: (Temporal graphs $|$ Streaming updates $|$ Multi-layer / Multiplex networks)
		\emph{Temporal graphs} track how community structures evolve over time using incremental optimization or temporal regularization.
		\emph{Streaming updates} process edges or nodes arriving in sequence without full recomputation.
		\emph{Multi-layer / Multiplex networks} detect communities across multiple interaction types or layers simultaneously, preserving inter-layer dependencies.	
\end{itemize}

\section{Preliminaries and Problem Variants}

A clear understanding of the community detection problem requires a precise formulation of the fundamental concepts involved. Although these notions are widely used in the literature, they are often introduced informally or with varying terminology. This section establishes the mathematical foundations that will be used throughout the paper, providing consistent definitions of networks, communities, and modularity, as well as a formal statement of the community detection problem.

\begin{definition}[Undirected network]
An undirected network is a graph $G = (V,E)$ where $V$ is the set of vertices and $E$ is the set of edges. Denoting $n=|V|$, the topological structure of $G$ can be defined by an $n \times n$ symmetric adjacency matrix $A = (a_{ij})_{n \times n}$, where $a_{ij} = 1$ if there is an edge between vertices $i$ and $j$ in $E$, and $a_{ij} = 0$ otherwise.
\end{definition}

\begin{definition}[Directed network]
A directed network is a directed graph $G = (V,E)$ where $V$ is the set of vertices and $E$ is the set of arcs. Denoting $n=|V|$, the topological structure of $G$ can be defined by an $n \times n$ adjacency matrix $A = (a_{ij})_{n \times n}$, where $a_{ij} = 1$ if $(i, j) \in E$ there is an arc from vertex i to vertex j in $E$, and $a_{ij} = 0$ otherwise.
\end{definition}

\begin{definition}[Weighted network]
A weighted network $(G, w)$ consists of a network $G = (V, E)$ and a function $w : E \to \R$ that assigns weights to the edges/arcs of $G$. The weight matrix $W = (w_{ij})_{n \times n}$ holds the weight of the edge/arc between $i$ and $j$, and by convention $w_{ij} = 0$ if there is no edge/arc between $i$ and $j$.
\end{definition}

\textbf{Remark:} An unweighted network $G$ can be thought of as a weighted network where the weight of every edge is the same.

\begin{definition}[Community in a network]
Given a network $G=(V,E)$, a community is a subset $C \subseteq V$ intended to represent a group of vertices that are more densely connected internally than with the rest of the network, according to a given quality function or structural criterion.
\end{definition}

\textbf{Remark:} There are at most $2^{|V|}$ different communities in a network $G = (V, E)$.
In general, a maximum number of communities $n_c \in \N: n_c \leq 2^{|V|}$ is imposed. 

\begin{definition}[Assignment of $V$ into communities]
	An assignment of $V$ into $n_c$ communities is a matrix $x \in \{0,1\}^{|V|\times n_c}$ (where $x_{ik} = 1 \Leftrightarrow \text{ vertex } i$ belongs to community $k\in\{1,...,n_c\}$) such that:
	\begin{enumerate}
		\item $\sum\limits_{k=1}^{n_c} x_{ik}\geq 1,  \, \forall i\in V$ (each vertex is assigned to at least one community).
		\item $(x_{ik})_{i\in V}\neq (x_{ik'})_{i\in V}\,\, \forall k,k'\in \{1,...,n_c\}: k\neq k' \wedge \sum\limits_{i\in V} x_{ik}\geq 1 \wedge \sum\limits_{i\in V} x_{ik'}\geq 1$ (any two distinct, non-empty communities must have different sets of vertices).
	\end{enumerate}

	In addition, the assignment is either one of the two following:
	\begin{itemize}
		\item  disjoint $\Leftrightarrow \sum\limits_{k=1}^{n_c} x_{ik}= 1,  \, \forall i\in V$ (each vertex is assigned to exactly one community).
		\item  overlapping $\Leftrightarrow \exists i\in V: \sum\limits_{k=1}^{n_c} x_{ik}> 1$ (at least one vertex is assigned to more than one community).
	\end{itemize}
\end{definition}

\begin{definition}[Fuzzy assignment of $V$ into communities] 
A fuzzy assignment of $V$ into $n_c$ communities for a given $\lambda \in [0, 1)$ (that denotes the threshold required for a vertex to belong to a community), 
is a pair $u, x$ where 
$u \in [0,1]^{|V|\times n_c}$ (represents the degree of membership of a vertex to a community), 
and $x \in \{0,1\}^{|V|\times n_c}$ is an assignment of $V$ into $n_c$ communities,
such that:
\begin{enumerate}
	\item $\sum\limits_{k=1}^{n_c} u_{ik}= 1,  \, \forall i\in V$ (the degree of membership of a vertex is normalized among all communities). 
	\item $x_{ik}=  \begin{cases} 
	1& \text{if } u_{ik} > \lambda \\
	0 & \text{if } u_{ik} \leq \lambda
\end{cases}
$ (a vertex belongs to a community $\Leftrightarrow$ its degree of membership is higher than $\lambda$). 
\end{enumerate}
\end{definition}

\textbf{Remark:} For any positive integer $p$, if $\dfrac{1}{p+1} \leq \lambda$, then each vertex can belong to at most $p$ communities. 
Thus, if $0.5 \leq \lambda $ the fuzzy assignment is disjoint.

\textbf{Example:} According to the definitions previously given for assignments of vertices into communities, we illustrate with an example the relationship between $x$, $u$ and $\lambda$. Given $V=\{1,...,5\}$ and $n_c=2$, Figure 1 shows a disjoint assignment of $V$ into communities. The associated matrix $x$ is linked to the matrix $u$ and the value $\lambda$:\\
$$x=\begin{pmatrix}
	1 & 0 \\
	1 & 0 \\
	1 & 0 \\
	0 & 1 \\
	0 & 1
\end{pmatrix}
\hspace{1cm}
u=\begin{pmatrix}
	\lambda+\mu_1 & 1- \lambda-\mu_1\\
	\lambda+\mu_2             & 1- \lambda-\mu_2 \\
	\lambda+\mu_3             & 1- \lambda-\mu_3 \\
	1- \lambda-\mu_4             & \lambda+\mu_4 \\
	1- \lambda-\mu_5             & \lambda+\mu_5
\end{pmatrix}
\text{ with } \lambda < 1 \text{ and } \mu_i\in (0,1-\lambda], \forall i\in V.$$
Analogously, Figure 2 shows an overlapping assignment of $V$ into communities where:
\begin{itemize}
    \item vertex $3$ is now also included in community 2: $x_{3,2} = 1$.
    \item the upper bound for $\lambda$ changes: $\lambda < \dfrac{1}{2}$. If this were not the case, then no vertex could belong to more than one community (see previous remark).
    \item the upper bound for $\mu_3$ changes: $\mu_3\in (0,1 - 2\lambda]$. If the degree of membership for vertex 3 to one community were too high, then it would not exceed the threshold $\lambda$ for the other community.
\end{itemize}

\begin{figure}[h]
    \begin{minipage}{0.55\textwidth}
        \centering
        \includegraphics[scale=1.35]{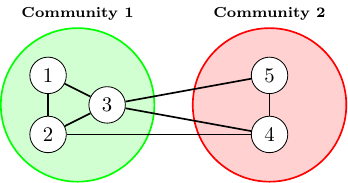}
        \caption{Example of disjoint communities.}
    \end{minipage}%
    \hspace{2mm}
    \begin{minipage}{0.45\textwidth}
        \centering
        \includegraphics[scale=1.35]{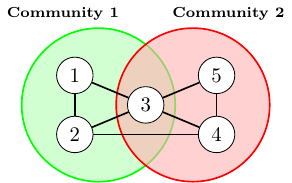}
        \caption{Example of overlapping communities.}
    \end{minipage}
\end{figure}

\begin{definition}[Modularity function]
	Let $(G,w)$ be a weighted network with $|V|=n$, let $n_c \in \mathbb{N}$, and 
	let $\mathcal{A}$ the feasible set of assignment of vertices into communities.
	
	A modularity (or quality) function is a mapping
	\[
	Q : (G,w,a) \mapsto \mathbb{R}
	\]
	that evaluates the quality of the assignment $a\in \mathcal{A}$. Larger values of $Q$ indicate better community structures, in the sense that vertices assigned to the same community are more strongly connected (or more similar) than expected under a suitable baseline or null model.
\end{definition}

\textbf{Example	[Newman--Girvan modularity]:}
One of the most widely used quality functions is the modularity proposed by \citet{Newman2004}. Given a partition of the vertex set into $n_c$ communities (i.e., a disjoint assignment), let $E$ be the $n_c \times n_c$ matrix whose entry $e_{kk'}$ represents the fraction of total edge weight linking vertices in community $k$ to vertices in community $k'$. Let $a_k = \sum_{k'=1}^{n_c} e_{kk'}$ denote the fraction of edge weight incident to community $k$.

The modularity is defined as:
\begin{equation}
	Q = \sum_{k=1}^{n_c} \left( e_{kk} - a_k^2 \right).
\end{equation}

This formulation is induced by a particular assignment matrix $X$, where each vertex belongs to exactly one community. In this sense, the quantities $e_{kk'}$ and $a_k$ are aggregated measures derived from the underlying assignment.

The function $Q$ compares the observed fraction of intra-community edges with the expected fraction under a random null model preserving vertex degrees.

\begin{definition}[Community detection problem]
	Let $(G,w)$ be a weighted network, let $n_c \in \mathbb{N}$, and let $Q(G,w,a)$ be a modularity (quality) function.
	
	The community detection problem consists of finding an assignment $a$ such that $Q$ is maximized, that is:
	\begin{equation}
		\max_{a\in \mathcal{A}} \; Q(G,w,a),
	\end{equation}
	where $\mathcal{A}$ denotes the feasible assignment space according to the chosen model (e.g., disjoint, overlapping).
	
	Additional constraints may be imposed on $\mathcal{A}$ depending on the application.
\end{definition}

\section{Literature Review on Modularity-Based Community Detection}

Modularity, introduced by \citet{Newman2004}, has become one of the most influential objective functions in community detection, serving both as a quality measure for evaluating partitions and as a direct optimization criterion for identifying community structure. The measure quantifies the deviation between observed intra-community connectivity and that expected under a degree-preserving random graph model. Its interpretability, analytical structure, and empirical effectiveness have established modularity maximization as a central problem in network science and related optimization domains.

From a mathematical standpoint, modularity maximization is a discrete, nonconvex optimization problem with a quadratic objective function and combinatorial assignment constraints. The combinatorial nature of the problem renders exact optimization computationally challenging; \citet{Aloise2010} even proved NP-hardness for the modularity function proposed in \citet{Newman2006}. These characteristics have stimulated two main lines of research. The first develops algorithmic strategies based on heuristics, relaxations, and spectral techniques, primarily aimed at scalability and computational efficiency. The second formulates the problem explicitly within mathematical programming frameworks, seeking exact or provably optimal solutions.

This section provides a structured review of the literature on modularity-based community detection. We distinguish between approaches designed for disjoint community structures and those addressing overlapping communities. Within each category, we further separate contributions that rely on mathematical programming formulations from those that adopt alternative algorithmic strategies. This organization highlights both the methodological diversity of modularity optimization and the comparatively limited number of contributions that formulate the problem within an explicit optimization framework. A comprehensive summary of the modularity-based algorithmic contributions in the reviewed literature, categorized by supported network characteristics and solution methods, is provided in Table \ref{tab:papers}.

\subsection{Disjoint Community Detection}

\subsubsection{Algorithmic Approaches}

The foundational contribution of \citet{Newman2004} established modularity as a quantitative objective for assessing network partitions and proposed a divisive procedure based on edge betweenness centrality. While the algorithm itself was not formulated as a direct modularity maximization method, the introduction of the modularity function naturally led to its interpretation as an optimization criterion.

Following the introduction of modularity as an optimization criterion, several greedy strategies were proposed to improve computational scalability. In \citet{Clauset2004}, a hierarchical agglomerative algorithm was introduced that iteratively merges communities to produce the largest increase in modularity at each step. The method significantly reduces computational complexity compared to earlier divisive procedures and enabled the analysis of substantially larger networks. Among these scalable heuristics, the Louvain method introduced in \citet{0DBlondel2008} has become one of the most widely used algorithms for modularity maximization. The method alternates between local node reassignment moves that increase modularity and aggregation steps that construct a reduced network whose vertices correspond to detected communities. This multi-level scheme yields high-quality partitions with near-linear empirical complexity on large networks. This paradigm was recently advanced by the Leiden algorithm \citep{Traag2019}, which refines the local moving phase to ensure that detected communities are internally well-connected, addressing structural instabilities occasionally found in Louvain partitions.

Spectral methods form a central class of approaches for modularity maximization. In \citet{White2005}, modularity optimization for weighted undirected graphs is addressed via spectral relaxation, resulting in two clustering algorithms. In \citet{Newman2006}, modularity was expressed in matrix form through the modularity matrix, casting modularity maximization as a quadratic optimization problem and enabling an eigenvector-based approximation algorithm. In \citet{Leicht2008}, the modularity function is generalized to unweighted directed networks, together with a corresponding spectral optimization procedure. These contributions establish a connection between modularity maximization and classical spectral partitioning techniques through relaxation of the discrete assignment constraints.

A further development concerns the introduction of a tunable resolution parameter. In \citet{Reichardt2006}, modularity is generalized by incorporating a parameter controlling the relative weight of the null model term, with community structure identified via simulated annealing. This extension allows the detection of communities at different scales and broadens the flexibility of modularity-based
optimization.

Beyond spectral and greedy techniques, metaheuristic approaches have been explored to search the combinatorial solution space more broadly. For instance, \citet{Tasgin2007} introduce a genetic algorithm designed for modularity maximization, and \citet{Nascimento2013} present a GRASP with path relinking algorithm for weighted graphs, combining randomized greedy construction with local search and a memory mechanism based on previously found high-quality solutions. Such methods aim to balance solution quality and computational effort without providing explicit optimality guarantees.

\citet{deSantiago2017Heur} address the scalability limitations of existing methods for modularity density maximization by introducing seven heuristics divided into three categories: a coarsening merger, moving node methods, and multilevel heuristics. The paper also establishes that the modularity density function cannot be used directly as a prioritizer in coarsening merger heuristics, as it provably merges cliques that should remain separate, and proposes an alternative prioritizer based on edge density.

Collectively, algorithmic approaches constitute the predominant methodology for modularity-based community detection, with emphasis on scalability and adaptability in large or structurally heterogeneous networks. A comprehensive survey of heuristic- and metaheuristic-based community detection algorithms is presented in \citet{Attea2021}. 

\subsubsection{Mathematical Programming Approaches}

In addition to spectral and heuristic procedures, modularity maximization has been modeled explicitly within mathematical programming frameworks that encode its combinatorial structure.

An integer linear programming formulation based on pairwise community membership variables is presented in \citet{Brandes2007}. Binary variables indicate whether two vertices belong to the same community, and linear transitivity constraints are enforced to ensure a valid partition of the network. Building on this, \citet{Agarwal2008} investigate modularity maximization through the lens of linear programming relaxations and rounding techniques. They demonstrate that while exact ILP provides a rigorous benchmark for small-to-medium networks, rounding-based heuristics derived from these relaxations can achieve near-optimal solutions where exact methods become computationally prohibitive.

Alternative formulations preserve the quadratic structure of the modularity objective. In \citet{Xu2007}, modularity maximization is expressed as a mixed-integer quadratic program (MIQP) using binary assignment variables to represent community membership directly. The quadratic objective captures the contribution of vertex pairs while maintaining explicit partition constraints.

Extension to directed graphs has further expanded the scope of mathematical programming in this domain. In \citet{Yang2016}, an integer linear programming (ILP) formulation is developed specifically for directed modularity. By adapting pairwise membership variables to the asymmetric directed null model, this framework provides exact benchmarks for directed networks. This allows for a rigorous evaluation of accompanying metaheuristics, such as simulated annealing, and highlights the increased complexity of the optimization landscape when edge directionality is preserved.

Decomposition-based strategies have also been investigated. In \citet{Aloise2010}, a column generation framework is proposed in which each column represents a candidate community and modularity is expressed as a sum of community-level contributions. The master problem selects a partition of the vertex set, while the pricing problem identifies communities with positive reduced cost. This approach leverages the additive structure of the modularity function and avoids the explicit enumeration of all possible partitions.

Nonlinear integer formulations have likewise been considered. In \citet{Bennett2012-2}, modularity maximization is modeled as a nonlinear integer program reflecting the quadratic dependence between vertex assignments. Similarly, in \citet{Xu2010}, a mixed-integer nonlinear programming (MINLP) model is introduced and solved through an iterative scheme based on multiple initial solutions. The method combines exact optimization modeling with multi-start refinement in order to improve practical solution quality.

\citet{Li2008} propose \textit{modularity density} (the $D$ value) as an alternative objective function that addresses the well-known \emph{resolution limit} of Newman and Girvan's modularity (the tendency of modularity maximization to overlook small communities by merging them into larger ones) and formulate its maximization as a nonlinear integer program. A generalized variant $D_\lambda$ is also introduced, where the parameter $\lambda$ interpolates between two known graph partition criteria---ratio association and ratio cut---enabling multi-resolution community analysis.

\citet{Cafieri2015} study the impact of adding cohesion conditions as explicit constraints to a mixed-integer quadratic programming formulation for modularity maximization. Five conditions of increasing strictness are considered --- strong, semi-strong, almost-strong, weak, and extra-weak --- and each is encoded as a set of linear constraints adjoined to the base model. Computational experiments on real-world networks show that stricter conditions substantially reduce modularity but can overcome the resolution limit, while the weak condition yields intuitive partitions with only moderate loss in objective value.

\citet{Costa2015} addresses the modularity density maximization problem, whose natural formulation is a binary nonlinear program due to the fractional structure of the objective. Several MILP reformulations are derived by linearizing bilinear terms via Fortet and McCormick inequalities \citep{Fortet1960, McCormick1976} and introducing auxiliary variables to handle the denominators, with a binary decomposition variant producing the most compact and computationally efficient model. Experiments on benchmark instances demonstrate that the best MILP reformulation solves instances orders of magnitude faster than nonlinear solvers applied to the original 0--1 NLP.

\citet{Costa2017} extend this line of work by providing complete MILP formulations that eliminate the auxiliary nonlinear subproblems still required in \citet{Costa2015} to compute bounding parameters. Exact MILP reformulations of those subproblems are derived using binary expansion of integer variables and dual reformulation of multilinear terms, and it is additionally proved that the continuous relaxations of the auxiliary problems always yield integer optimal solutions. The proposed formulations reduce resolution time by up to two orders of magnitude compared to the previous approach.

\citet{Sato2019} present an exact algorithm for the modularity density maximization problem on undirected graphs. They formulate the problem as a set-partitioning integer linear program and solve it via a branch-and-price framework, in which the column generation subproblem is expressed as a MILP --- a simpler formulation than the integer quadratic programming approach proposed by \citet{deSantiago2017Ex}. Two acceleration techniques are incorporated: a set-packing relaxation, which replaces the set-partitioning constraint with a weaker set-packing constraint for most vertices and restores it selectively, and a multiple-cutting-planes-at-a-time strategy, which generates several non-overlapping columns per column generation iteration. Together, these techniques allow the authors to solve benchmark instances with over 100 vertices to proven optimality in around four minutes on a standard PC, finding new optimal solutions for several instances that had remained unsolved.

Taken together, these contributions illustrate how mathematical programming can bring rigorous exactness to modularity-based community detection, through a variety of formulations --- linear, quadratic, and nonlinear --- and solution strategies such as column generation and branch-and-price. Beyond providing optimal solutions for benchmark instances, these frameworks serve as a foundation for evaluating heuristic methods and for extending modularity optimization to richer problem settings such as directed graphs and alternative objective functions like modularity density.

\subsubsection{Structural Properties of Modularity}

\citet{Fortunato2007} showed that modularity maximization may aggregate small communities within larger structures depending on global network characteristics, a behavior known as the resolution limit and a direct consequence of the global null model embedded in the modularity definition. \citet{Good2010} identify a complementary issue known as \emph{degeneracy}: the modularity function typically admits an exponential number of structurally diverse partitions with near-optimal objective value and no clear global maximum. As a consequence, different algorithms applied to the same network may return very different partitions with similarly high modularity scores, making it difficult to determine which partition, if any, reflects the true community structure.

The identification of these phenomena motivated subsequent work on refined optimization approaches, such as mathematical programming formulations that incorporate iterative or multi-start strategies to improve detection across different structural scales \citep{Xu2010}.

Thus, the resolution limit and degeneracy should be understood as structural characteristics of modularity that inform model refinement and the cautious interpretation of results, rather than as deficiencies that invalidate its use.

\subsection{Overlapping Community Detection}

The original modularity formulation is defined over partitions and therefore assumes disjoint communities. However, the foundational work of \citet{Palla2005} via clique percolation highlighted that real-world networks actually exhibit interwoven, overlapping structures. While their specific approach was based on local dense substructures rather than modularity, it catalyzed a broader shift in how communities are modeled. Driven by this new paradigm, researchers recognized a critical gap in optimization models: as observed by \citet{Chen2010}, classical modularity simply could not accommodate these multiple memberships. This realization directly motivated the development of the extended formulations capable of modeling overlap.

Compared with the disjoint setting, the literature on overlapping modularity maximization is more limited, particularly within a mathematical programming framework. This gap suggests opportunities for further methodological development. In the following, we review the existing contributions by distinguishing between algorithmic approaches and mathematical programming models.

\subsubsection{Algorithmic Approaches}

Early approaches to overlapping modularity maximization relied primarily on heuristic or clustering-based techniques. Among these early algorithmic methods, the spin-glass framework of \citet{Reichardt2006} provided a notable foundation for relaxing strict boundaries. As discussed previously, their model optimizes a generalized modularity via the heuristic of simulated annealing; however, a crucial additional feature of this physical mapping is that it inherently accommodates overlapping structures by allowing nodes to exist in mixed spin states.

Building upon this concept of flexible cluster assignments, subsequent methods explicitly modeled continuous overlapping memberships. In \citet{Zhang2007}, the authors consider weighted graphs with nonnegative edge weights and model overlapping communities through a fuzzy partition of the vertex set into $k$ groups. Instead of assigning each vertex to a single community, membership degrees are introduced, allowing vertices to belong to multiple communities simultaneously. The resulting objective function generalizes modularity to this fuzzy setting. To compute a solution, the authors combine the modularity concept with spectral relaxation \citep{White2005} and the fuzzy $k$-means clustering method \citep{Dunn1973}, leading to a global heuristic algorithm.

To overcome the computational demands of global matrix relaxations and fuzzy partitions, other algorithmic developments shifted toward localized heuristics. Returning to \citet{Chen2010}, the authors proposed a local expansion algorithm for weighted networks. Instead of computing a global partition, their approach identifies crisp overlapping communities from the bottom up by initializing a partial community at a high-strength node and iteratively adding adjacent nodes based on a local belonging degree metric.

Another extension of modularity to overlapping communities is proposed in \citet{Nicosia2009}. Starting from the modularity definition for directed networks introduced by \citet{Leicht2008}, the authors note that the classical formulation assumes sharply separated communities and therefore does not account for overlaps. They introduce a generalized modularity function for directed graphs with overlapping communities and propose a genetic algorithm to maximize the resulting objective.

Building further on the use of evolutionary algorithms to optimize generalized modularity, \citet{Tsung2020} addressed a critical limitation in standard fuzzy formulations, noting that strict fuzzy partitions where a node's memberships must sum exactly to one often fail to provide a strong enough signal to uncover overlapping structures. To resolve this, they reformulated the task as a node weight allocation problem. Their approach introduces a vertex-reweighting mechanism that permits the total membership weight of a shared node across all its communities to exceed one, thereby amplifying the structural presence of overlaps within the objective function. To optimize this continuous allocation, the authors employ a customized genetic algorithm paired with three localized refinement strategies—thresholding negligible memberships, merging heavily overlapping clusters, and proportionally reweighting nodes—to fine-tune the final community boundaries.

A broader survey of heuristic and metaheuristic approaches to overlapping community detection, among other settings, is provided in \citet{Attea2021}.

\subsubsection{Mathematical Programming Approaches}

Optimization-based formulations for overlapping modularity maximization are comparatively scarce compared to algorithmic ones. In \citet{Bennett2012-2}, the authors propose a nonlinear integer programming formulation for modularity maximization that captures the quadratic dependence between vertex assignments. The resulting MINLP provides an exact optimization framework for community detection and can be solved using standard mathematical programming solvers.

More recently, \citet{Benati2022} propose an integer linear programming formulation for maximizing the fuzzy modularity function introduced in \citet{Zhang2007}. The resulting model, denoted F-MOD, provides an exact optimization framework for the problem. The authors further observe that the formulation may allow communities to be fully contained within others. To address this issue, they introduce an alternative model, F-MOD-NI, that prevents such nested community structures. In addition, they propose an alternative modularity definition for weighted graphs to prevent the creation of many large communities that contain almost all vertices. It is also based on a different random graph interpretation: instead of a single random graph, the model considers $n_c$ random graphs whose edge weights correspond to averages of the membership degrees to the communities. A corresponding integer linear programming formulation is derived for this variant as well. The resulting models are solved using mathematical programming techniques, while heuristic procedures are employed for larger instances.

\begin{table*}[h!]
\centering
\caption{Summary of selected community detection literature, classified by supported network structure (D=Disjoint, O=Overlapping, We.=Weighted, Di.=Directed), mathematical programming formulation (Math. Prog.), and solution method.}
\label{tab:papers}
\small
\setlength{\tabcolsep}{6pt}
\begin{tabular}{c c c c c C{5.4cm}}
\toprule
Paper                                 & \makecell{D/O} & We. & Di. & \makecell{Math. Prog.} & \makecell{Solution method} \\
\midrule
\cite{Newman2004} & D & \Cross          & \Cross   & \Cross           & Girvan-Newman heuristic (divisive) \\
\cite{Clauset2004}       & D                       & \Cross          & \Cross   & \Cross           & CNM heuristic (greedy agglomeration)                \\
\cite{0DBlondel2008} & D & \Checkmark  & \Cross & \Cross  & Louvain heuristic (local moving)  \\
\cite{Traag2019}         & D                       & \Checkmark      & \Checkmark & \Cross           & Leiden heuristic (local moving with refinement) \\
\cite{White2005}         & D                       & \Checkmark          & \Cross   & \Cross          & spectral clustering \\
\cite{Newman2006}        & D                       & \Cross          & \Cross   & \Cross           & spectral optimization \\  
\cite{Leicht2008}         & D                       & \Cross          & \Checkmark   & \Cross          & spectral optimization \\
\cite{Reichardt2006}     & D, O       & \Checkmark & \Checkmark & \Cross           & simulated annealing \\
\cite{Tasgin2007} & D                       & \Cross          & \Cross    & \Cross         & genetic algorithm         \\
\cite{Nascimento2013}    & D                       & \Cross          & \Cross     & \Cross           & GRASP with path relinking \\
\cite{deSantiago2017Heur}    & D & \Checkmark      & \Cross     & \Cross & LNM/MCN heuristics (local moving and coarsening) \\
\cite{Brandes2007} & D                       & \Cross          & \Cross    & \Checkmark         & ILP, greedy          \\
\cite{Agarwal2008} & D & \Cross     & \Cross     & \Checkmark & LP rounding \\ 
\cite{Xu2007}         & D                       & \Cross          & \Cross  & \Checkmark           & MIQP                 \\
\cite{Yang2016}                                        &  D                       &  \Cross           & \Checkmark              & \Checkmark & MINLP + MIP                     \\
\cite{Aloise2010} &  D  & \Cross \Checkmark   & \Cross \Checkmark  & \Checkmark & column generation   \\
\cite{Bennett2012-2} &  D, O  &  \Checkmark & \Cross & \Checkmark  &  MINLP                     \\
\cite{Xu2010} & D &  \Cross &   \Cross & \Checkmark &  MIQP + iterative improv. \\
\cite{Li2008}            & D    & \Cross     & \Cross & \Checkmark & INLP / kernel $k$-means    \\ 
\cite{Cafieri2015}       & D & \Cross & \Cross & \Checkmark & MIQP + cohesion constraints \\
\cite{Costa2015}         & D & \Cross & \Cross & \Checkmark & MILP              \\
\cite{Costa2017}         & D & \Cross & \Cross & \Checkmark & MILP \\
\cite{Sato2019} & D & \Cross & \Cross & \Checkmark & branch-and-price \\
\cite{deSantiago2017Ex} & D & \Cross & \Cross & \Checkmark & column generation \\
\cite{Zhang2007} & O & \Checkmark  & \Cross & \Cross & fuzzy $k$-means clustering \\
\cite{Chen2010} & O & \Checkmark & \Cross & \Cross & local expansion heuristic \\
\cite{Nicosia2009} & O & \Cross & \Checkmark & \Cross & genetic algorithm \\
\cite{Tsung2020} & O & \Cross & \Cross & \Cross & genetic algorithm with refinement strategies \\
\cite{Benati2022} &  O &  \Cross \Checkmark & \Cross & \Checkmark & MILP, heuristics \\
\bottomrule
\end{tabular}
\end{table*}

\section{Modularity Functions}

Modularity is one of the most widely used objective functions for community detection, measuring the extent to which the observed connectivity within communities differs from that expected under a suitable null model. Over time, several formulations of modularity have been proposed to accommodate different types of networks and community structures.

In this section, we review a number of modularity functions used in the literature presenting their mathematical expressions and briefly discussing their underlying modeling assumptions. The starting point for most modularity formulations is the definition introduced by \citet{Newman2004}, which we present below. 

\subsection{Disjoint Communities}

The classical formulation of modularity assumes that the vertex set is partitioned into disjoint communities, so that each vertex belongs to exactly one group. Under this assumption, community detection amounts to identifying a partition that maximizes the difference between the observed number of intra-community edges and the expected number under a suitable null model.

The functions reviewed in this subsection adopt this hard partition framework and differ mainly in how the null model is defined or extended to account for weighted or directed networks. In the following, we present several formulations of modularity for disjoint communities, focusing on their mathematical expressions and underlying modeling assumptions.

\subsubsection{\citet{Newman2004}}

\citet{Newman2004} introduce modularity as a quality measure to evaluate a graph's partition into disjoint communities, assuming an undirected and unweighted network. The central idea is to compare the observed number of edges within communities with the expected number under a random graph model that preserves the degree distribution of the original network. In this null model, edges are placed at random subject to the constraint that each vertex retains its degree, so that the expected connectivity reflects only the degree sequence and not any community structure.

Define the matrix $E \in \mathbb{R}^{n_c \times n_c}$, where each entry $e_{kk'}$ represents the proportion of edges in the graph that connect vertices in community $k$ to vertices in community $k'$. Let $a_k = \sum_{k'=1}^{n_c} e_{kk'}$ denote the fraction of edges that have at least one endpoint in community $k$.

The modularity function is then defined as
\begin{equation}
\label{eq:modularity:newman}
Q = \sum_{k=1}^{n_c} \left( e_{kk} - a_k^2 \right).
\end{equation}
The term $e_{kk}$ measures the observed fraction of edges within community $k$, while $a_k^2$ corresponds to the expected fraction of such edges under the null model. Thus, modularity quantifies the extent to which the observed intra-community connectivity exceeds what would be expected at random, with higher values indicating stronger community structure.

\subsubsection{\citet{White2005}}

\citet{White2005} consider undirected graphs with nonnegative edge weights represented by a symmetric matrix $W$. Given a partition $P_{n_c}$ of the vertex set into $n_c$ communities, the modularity function is defined as
\begin{equation}
    \sum_{k=1}^{n_c} \left[
        \frac{W(V_k, V_k)}{W(V, V)} -
        \left( \frac{W(V_k, V)}{W(V, V)} \right)^2
    \right],
    \label{eq:modularity:white}
\end{equation}
where $W(V', V'') = \sum_{i \in V',\, j \in V''} w_{ij}$ denotes the total weight of edges between subsets $V'$ and $V''$.

The term $\frac{W(V_k, V_k)}{W(V, V)}$ represents the fraction of total edge weight that lies within community $k$, while $\frac{W(V_k, V)}{W(V, V)}$ is the fraction of edge weight incident to vertices in $k$. As in the original formulation of \citet{Newman2004}, the second term corresponds to the expected fraction of intra-community weight under a null model in which edge endpoints are placed independently according to their weighted degrees. Thus, the modularity measures the deviation between observed and expected intra-community connectivity in the weighted setting.

\subsubsection{\citet{Newman2006}}

\citet{Newman2006} also introduce an alternative modularity function, which is linear rather than quadratic:
\begin{equation}
\label{Newman-lineal}
\dfrac{1}{2m} \sum\limits_{i,j \in V} \left( a_{ij} - \dfrac{d_i d_j}{2m} \right) \delta(k_i, k_j)
\end{equation}
where $m$ is the number of edges, $a_{ij} = 1$ if there is an edge between vertices $i$ and $j$, and $a_{ij} = 0$ otherwise, $d_v$ is the degree of $v$ in $G$, and $\delta(k_i, k_j) = 1$ if $i$ and $j$ belong to the same community and $0$ otherwise.

\subsubsection{\citet{Leicht2008}}

\citet{Leicht2008} extend the Newman--Girvan modularity to directed networks by adapting the null model to account for edge directionality. In this setting, the adjacency matrix $A$ is generally asymmetric, and each node $i$ is characterized by its in-degree $d_i^{in}$ and out-degree $d_i^{out}$. The modularity function is defined as

\begin{equation}
    \frac{1}{m} \sum_{i,j \in V}
    \left(
        a_{ij} - \frac{d_i^{out} d_j^{in}}{m}
    \right)
    \delta(k_i, k_j).
    \label{eq:modularity:leicht}
\end{equation}

The key idea is to replace the undirected null model $\frac{d_i d_j}{2m}$ by a directed counterpart $\frac{d_i^{out} d_j^{in}}{m}$, which represents the expected number of edges from node $i$ to node $j$ in a random directed graph that preserves both in-degree and out-degree sequences. This modification ensures that the modularity measure properly accounts for the asymmetry of directed networks while retaining the interpretation of comparing observed connectivity with a degree-preserving null model.

\subsubsection{Modularity Density \citep{Li2008, Costa2015}}

\citet{Li2008} introduce the modularity density function as an alternative measure for evaluating the quality of a partition, motivated in part by the resolution limit associated with classical modularity. \citet{Costa2015} later proposes a more compact equivalent formulation that is more amenable to mathematical programming approaches. Given a partition of the vertex set into $n_c$ communities $V_1, \dots, V_{n_c}$, the modularity density is defined as
\begin{equation}
    \sum_{\substack{k=1 \\ |V_k| > 0}}^{n_c} \frac{2 m_k - \overline{m}_k}{|V_k|},
    \label{eq:modularity_density}
\end{equation}
where $m_k$ denotes the number of edges with both endpoints in $V_k$, $\overline{m}_k$ is the number of edges with exactly one endpoint in $V_k$, and $|V_k|$ is the number of vertices in community $k$.

The quantity $2m_k - \overline{m}_k$ captures the balance between internal connectivity and external connections of community $k$, while the normalization by $|V_k|$ accounts for the size of the community. In this way, modularity density evaluates the cohesion of each community at a local scale, rather than aggregating contributions globally as in standard modularity.

\subsection{Overlapping Communities}

The functions reviewed here relax the hard partition assumption, allowing nodes to belong to multiple communities simultaneously. The node set of community $k$ is then defined as $V_k = \{i \in V \mid u_{ik} > \lambda\}$ for some threshold $\lambda > 0$.

Among the formulations for overlapping communities, we begin with the fuzzy modularity proposed by \citet{Zhang2007}.

\subsubsection{\citet{Zhang2007}}
\citet{Zhang2007} extend the Newman--Girvan modularity to overlapping communities by incorporating fuzzy membership coefficients $u_{ik} \in [0,1]$, which capture the degree to which node $i$ belongs to community $k$. The network is represented by a non-negative symmetric weight matrix $W$, where $w_{ij} \geq 0$ denotes the weight of the edge between nodes $i$ and $j$, with $w_{ij} = 0$ indicating the absence of an edge. Define the following weighted interaction measures:
\begin{align}
    U(V_k, V_k) &= \sum_{i,\, j \in V_k}
        \frac{u_{ik} + u_{jk}}{2}\, w_{ij}, \label{eq:zhang:acc}\\[4pt]
    U(V_k, V) &= U(V_k, V_k) +
        \sum_{\substack{i \in V_k \\ j \notin V_k}}
        \frac{u_{ik} + (1 - u_{jk})}{2}\, w_{ij}, \label{eq:zhang:acv}
\end{align}

where all sums range over ordered pairs. The modularity is then
\begin{equation}
    \sum_{k=1}^{n_c} \left[
        \frac{U(V_k, V_k)}{W(V, V)}
        - \left( \frac{U(V_k, V)}{W(V, V)} \right)^{\!2}
    \right].
    \label{eq:modularity:zhang}
\end{equation}

\subsubsection{\citet{Nicosia2009}}

\citet{Nicosia2009} extend the modularity for directed networks of \citet{Leicht2008} to overlapping communities. The contribution of each edge $(i,j)$ to the modularity of community $k$ is weighted by $\beta_{l(i,j),k} = \varphi(u_{ik}, u_{jk})$, where $\varphi$ is a function left deliberately general; possible choices include the product, the maximum, or the average of the two belonging coefficients. The null model is then defined by computing, for each node, the expected belonging coefficient of its outgoing and incoming edges to community $k$:
\begin{align}
    \beta^{out}_{l(i,*),k} &= \frac{\sum_{j \in V}
        \varphi(u_{ik},\, u_{jk})}{|V|}, \label{eq:nicosia:bout}\\[4pt]
    \beta^{in}_{l(*,j),k} &= \frac{\sum_{i \in V}
        \varphi(u_{ik},\, u_{jk})}{|V|}. \label{eq:nicosia:bin}
\end{align}

If $\beta_{l(i,j),k} = \varphi(u_{ik}, u_{jk})$ is the weight of the contribution of $l(i,j)$ to the modularity of community $k$, then the modularity function is
\begin{equation}
    \frac{1}{m} \sum_{k=1}^{n_c} \sum_{i,j \in V}
    \left(
        \beta_{l(i,j),k}\, a_{ij}
        - \frac{\beta^{out}_{l(i,*),k}\, d^{out}_i\,
                \beta^{in}_{l(*,j),k}\, d^{in}_j}{m}
    \right).
    \label{eq:modularity:nicosia}
\end{equation}

\subsubsection{\citet{Benati2022}}

\citet{Benati2022} build on the formulation of \citet{Zhang2007}, but observe that such formulation is not appropriate for overlapping community detection in practice. The cross-boundary terms in $A(V_k, V)$ involve the complement $1 - u_{jk}$ for nodes $j \notin V_k$; since these terms appear subtracted in~\eqref{eq:modularity:zhang}, they introduce negative contributions that depend on non-membership. As a consequence, nodes with low membership in a community still generate large penalties, which creates an imbalance in the objective and encourages the formation of excessively large communities, leading optimal solutions to consist of communities containing almost all vertices. \citet{Benati2022} address this by removing those negative terms, retaining only the positive contributions. Let $W = \sum_{(i,j) \in E} w_{ij}$ denote the total edge weight and $W_v = \sum_{l \in N(v)} w_{vl}$ the weighted degree of node $v$. The resulting modularity function is
\begin{equation}
    \frac{1}{2W} \sum_{k=1}^{n_c} \sum_{i,j \in V_k}
    \left( w_{ij} - \frac{W_i W_j}{2W} \right)
    \frac{u_{ik} + u_{jk}}{2},
    \label{eq:modularity:benati}
\end{equation}

\section{Datasets Used in Community Detection Research}

The study of community detection is strongly guided by the networks to which these methods are applied. Datasets provide not only benchmarks for evaluating algorithmic performance but also a window into the structural characteristics and challenges inherent in the community detection problem. By examining a range of networks, from synthetic graphs with controlled properties to real-world systems drawn from social, biological, and technological domains, researchers can better understand the diversity of community structures and the requirements that effective detection methods must satisfy.

Different datasets highlight different aspects of the problem. Synthetic networks are typically designed with a known ground truth, which allows for precise evaluation of an algorithm's ability to recover communities under varying conditions such as network size, density, community size distribution, and mixing parameters. They serve as a testbed to study the robustness, scalability, and sensitivity of methods to structural variations. In contrast, real-world networks illustrate practical challenges that arise in community detection, including heterogeneous degree distributions, overlapping communities, noise in edge formation, and the absence of complete ground-truth labels. These networks motivate the development of methods that are not only accurate but also interpretable and adaptable to complex, irregular structures.

In this work, we treat synthetic benchmark networks separately from real-world datasets, which are organized according to their primary application domains. The following subsections introduce the domains considered and summarize representative datasets within each category.

The application domains considered in this survey are not completely mutually exclusive, and the assignment of a dataset to a particular domain is not always unique. Several networks can naturally be interpreted from multiple perspectives depending on the modeling objective and the semantics emphasized, for example social interaction, information flow, or technological infrastructure. In such cases, we assign each dataset to a single domain for consistency and clarity of presentation, while acknowledging that alternative classifications may also be reasonable.

\textbf{Social Systems (\Soc):} Social networks constitute one of the most established benchmark domains for community detection, as they naturally exhibit densely connected groups reflecting social organization, affiliation, or interaction patterns. A classical example is Zachary's Karate Club Friendship Network \citep{0DZachary1977}, a small-scale real-world network of 34 individuals whose known split into two factions provides a widely used ground truth for validation. Another widely used benchmark is the American College Football Network \citep{Girvan2002}, where nodes represent teams, edges correspond to games played, and communities align closely with conference memberships. At a much larger scale, the Belgian Mobile Phone Network \citep{0DBlondel2008} represents millions of users connected through call interactions, illustrating the challenges of community detection in massive real-world communication graphs where external metadata (e.g., language regions) are often used as proxy ground truth. Other datasets considered in this domain include the Dolphin Social Network \citep{0DLusseau2003}, Jazz Musicians Collaboration Network \citep{0DGleiser2003}, Highland Tribes Network \citep{0DRead1954}, Enron Email Network \citep{0DKlimt2004}, Krebs' Political Books Network \citep{Newman2006}, Les Misérables Character Network \citep{0DKnuth1993}, Facebook Network \citep{0DMcAuley2012} and Twitter Network \citep{0DDeDomenico2013}.

\textbf{Biological, Ecological, and Healthcare Systems (\BioEco):} This combined category encompasses networks that capture interactions among biological entities, species within natural ecosystems, or similarities derived from clinical data. In these networks, community structure often corresponds to functional cellular modules, trophic organization, or specific patient subgroups. At the molecular level, the Protein-Protein Interaction Network \citep{0DDunn2005} models physical interactions between proteins to identify functional complexes. At the macroscopic ecological level, datasets like the Food Web of Marine Organisms \citep{0DBaird1989} encode predator-prey relationships to study modularity, stability, and energy flow. Additionally, animal behavioral networks, such as the Macaque Network \citep{0DSade1972} and the Zebra Communication Network \citep{0DSundaresan2007}, bridge behavioral and population dynamics. Finally, in the biomedical context, the Human Disease Network \citep{0DGoh2007} links disorders to their associated disease genes, where communities strongly align with 22 distinct physiological system pathologies.

\textbf{Communication and Infrastructure (\Comm):} Communication and technological networks describe infrastructure-level connectivity and information exchange systems, often characterized by large scale, sparsity, and heterogeneous degree distributions. The P2P Network \citep{0DRipeanu2002} models peer-to-peer file sharing systems, where communities can reflect clusters of peers with similar connectivity patterns. Another representative dataset is the US Airports Network \citep{0DOpsahl2010}, which models airports as nodes and flight connections as edges and exhibits pronounced hub structure and spatially induced clustering. Other datasets in this domain include Mobile Ad Hoc Networks \citep{0DHui2011} and Web Pages Network \citep{0DBlondel2008}.

\textbf{E-Commerce and Economics (\EComm):} E-commerce and economic networks represent interactions between users, products, or transactions, where communities often correspond to consumer segments, product categories, or market structures. The Amazon Products Network \citep{0DLeskovec2007} connects products that are frequently co-purchased, forming large graphs with implicit category-based community structure. The Epinions Trust Network \citep{0DRichardson2003} models trust relationships among users of a product review platform and is commonly used to analyze social influence and recommendation dynamics. Other datasets in this domain include the eBay Network \citep{0DReichardt2007} and the Retail Transaction Network from Chile \citep{0DRios2014}.

\textbf{Academia and Scientometrics (\Aca):} Academic and scientometric networks describe collaboration, citation, and communication patterns in scientific production, where communities often reflect research fields, institutional structures, or collaboration circles. A representative benchmark is the Santa Fe Institute Collaboration Network \citep{Girvan2002}, which models co-authorship relationships among scientists affiliated with the Santa Fe Institute and captures the interdisciplinary structure of research collaborations within a well-defined scientific community. Its moderate size and clear semantics make it a frequently used benchmark for validating community detection methods on collaboration networks. At a larger scale, the Microsoft Academic Graph \citep{0DWang2016} provides a massive heterogeneous network of papers, authors, venues, and citations, enabling large-scale analysis of scholarly structure and evolution. Additional datasets in this domain include the Physics E-print Archive Network \citep{0DNewman2001}, the Netscience Co-Authorship Network \citep{Clauset2004}, the Scientific Communication Network \citep{0DLeydesdorff2009} and the Cora Network \citep{0DSen2008}.

\textbf{Information and Semantics (\Info):} Information and semantic networks encode relationships between documents, concepts, or words, where communities often correspond to topical or semantic clusters. The Word Association Network of Cognitive Sciences \citep{0DNelson2004} represents associations between words collected from human experiments and is used to study semantic organization and cognitive structure through clustering. The Political Blogs Network \citep{0DAdamic2005} is a directed network of hyperlinks between weblogs on US politics recorded around the 2004 presidential election, where nodes are classified according to political leaning into two groups, providing a ground truth that reflects ideological alignment. Other datasets in this domain include the Stanford Web Graph and Google Web Graph \citep{0DLeskovec2009}, which model large-scale hyperlink structure among web pages.

We summarize commonly used datasets in Table \ref{tab:datasets}. Synthetic benchmark generators are excluded from the table, since they do not correspond to fixed datasets but rather to parametric families of graphs whose characteristics vary across experimental settings. For each real-world dataset, we provide its common abbreviation and full name, specify its application domain, and report known ground-truth communities when available. When ground truth is available, we also indicate the number of distinct labels (e.g., ``Yes: 2''). Additionally, we include descriptive statistics such as the number of vertices and edges, whether the network is directed or undirected, the year of release, and the accessibility of the data.

\begin{table*}[h]
\centering
\begin{minipage}{\textwidth}
\centering
\caption{Real-world datasets and their main structural characteristics.}
\label{tab:datasets}
\small
\setlength{\tabcolsep}{6pt}
\begin{tabular}{c c c c c c}
\toprule
\textbf{Dataset} 
& \makecell{\textbf{Application}\\\textbf{Domain}} 
& \textbf{Ground Truth} 
& $\lvert V \rvert$ 
& $\lvert E \rvert$ 
& \textbf{Directed} \\
\midrule

Zachary's Karate Club\tabcite{0DZachary1977} 
& \Soc 
& Yes: 2 
& 34 
& 78 
& No \\

American College Football\tabcite{Girvan2002}
& \Soc 
& Yes: 12 
& 115 
& 613 
& No \\

Belgian Mobile Phones\tabcite{0DBlondel2008} 
& \Soc 
& Yes: 2
& ~2.6M
& ~6.3M
& No \\

Dolphins\tabcite{0DLusseau2003} 
& \Soc 
& Yes: 2 
& 62 
& 159 
& No \\

Jazz Musicians\tabcite{0DGleiser2003} 
& \Soc 
& No 
& 198 
& 2742 
& No \\

Highland Tribes\tabcite{0DRead1954} 
& \Soc 
& Yes: 3
& 16
& 58
& No \\

Enron Emails\tabcite{0DKlimt2004}
& \Soc 
& No
& 36692
& 183831
& Yes \\

Krebs' Political Books\tabcite{Newman2006} 
& \Soc 
& Yes: 3
& 105
& 441 
& No \\

Les Misérables Character\tabcite{0DKnuth1993} 
& \Soc 
& No
& 77
& 254 
& No \\

Facebook\tabcite{0DMcAuley2012} 
& \Soc 
& Yes: 193
& 4309
& 88234 
& No \\

Twitter\tabcite{0DDeDomenico2013} 
& \Soc 
& No
& 456626
& ~14.8M
& Yes \\

Protein-Protein Interactions\tabcite{0DDunn2005} 
& \BioEco 
& No
& 2444
& 6271
& No \\

Food Web of Marine Organisms\tabcite{0DBaird1989} 
& \BioEco
& No
& 128
& 2106 
& Yes \\

Macaques\tabcite{0DSade1972} 
& \BioEco
& Yes: 3
& 16
& 111
& Yes \\

Zebras\tabcite{0DSundaresan2007} 
& \BioEco
& No
& 27
& 111
& No \\

Human Diseases\tabcite{0DGoh2007} 
& \BioEco
& Yes: 22
& 1419 
& 2738 
& No \\

Gnutella P2P\tabcite{0DRipeanu2002} 
& \Comm
& No
& 6301
& 20777 
& Yes \\

US Airports\tabcite{0DOpsahl2010} 
& \Comm
& No
& 1574 
& 28236 
& No \\

Mobile Ad Hoc\tabcite{0DHui2011} 
& \Comm
& Yes: 11
& 98
& 4226
& No \\

Web Pages\tabcite{0DBlondel2008} 
& \Comm
& No
& 39459
& 442247
& Yes \\

Amazon Products\tabcite{0DLeskovec2007} 
& \EComm
& Yes: 4
& 262111
& 1234877
& Yes \\

Epinions Trust\tabcite{0DRichardson2003} 
& \EComm
& No
& 75879
& 508837
& Yes \\

eBay\tabcite{0DReichardt2007} 
& \EComm
& No
& 11885 
& 22865
& No \\

Retail Transaction\tabcite{0DRios2014} 
& \EComm
& Yes: 18
& 11803 
& 148155
& No \\

Santa Fe Institute\tabcite{Girvan2002} 
& \Aca
& Yes: 7
& 118
& 197
& No \\

Microsoft Academic\tabcite{0DWang2016} 
& \Aca
& Yes: 19
& 2827700
& 48124312 
& Yes \\

Physics E-print Archive\tabcite{0DNewman2001} 
& \Aca
& No
& 34546 
& 421578
& Yes \\

Netscience Co-authorship\tabcite{Clauset2004}
& \Aca
& No
& 1589 
& 2742 
& No \\

Scientific Communication\tabcite{0DLeydesdorff2009}
& \Aca
& Yes: 14
& 6128
& 50716
& Yes \\

Cora\tabcite{0DSen2008}
& \Aca
& Yes: 7
& 2708
& 5429
& Yes \\

Word Association\tabcite{0DNelson2004}
& \Info
& No
& 10617 
& 72168 
& Yes \\

Political Blogs\tabcite{0DAdamic2005}
& \Info
& Yes: 2
& 1490
& 19090
& Yes \\

Stanford\tabcite{0DLeskovec2009}
& \Info
& No
& 281903 
& 2312497
& Yes \\

Google\tabcite{0DLeskovec2009}
& \Info
& No
& 875713 
& 5105039 
& Yes \\

\bottomrule
\end{tabular}
\end{minipage}
\end{table*}

\section{Discussion, Open Issues and Future Research Directions}

This paper provides a structured and unified perspective on the Community Detection Problem by integrating four complementary components: a comprehensive literature review, a multi-dimensional taxonomy of classification tags, a formal mathematical framework, and a consolidated view on evaluation through modularity functions and benchmark datasets.

First, the literature review highlights the progressive diversification of the field. Early works focused primarily on algorithmic strategies and structural definitions, while more recent contributions emphasize specialized settings such as multilayer, temporal, attributed, or learning-based networks. This evolution has led to a fragmented landscape in which different works adopt heterogeneous terminologies, modeling assumptions, and evaluation criteria.

Second, to address this fragmentation, we have proposed a unified taxonomy organized along six main dimensions: network structure, type of community, objective function, methodological approach, evaluation and validation, and application domain. This taxonomy consolidates the various classification criteria found in the literature into a coherent framework, enabling systematic comparison across methods and clarifying conceptual relationships that are often treated independently.

Third, we have introduced a general mathematical formalization of the problem based on assignment functions, allowing both crisp and fuzzy community memberships to be represented within a single framework. This formulation provides a unifying viewpoint in which disjoint, overlapping, and fuzzy community detection can be interpreted as particular cases of a broader assignment model. Furthermore, by defining modularity in a generalized way as a quality function over assignments, the framework accommodates the wide variety of objective functions proposed in the literature.

Fourth, we have reviewed and analyzed representative modularity functions, highlighting their underlying assumptions, implicit null models, and known limitations. This analysis reinforces the idea that the choice of quality function is not neutral, but rather determines the type of communities that can be detected. In particular, issues such as the resolution limit, degeneracy of optimal solutions, and sensitivity to network structure illustrate the need for careful selection or combination of objective functions.

Finally, we have examined commonly used benchmark datasets and evaluation practices in the literature. This review reveals a lack of standardization in experimental settings, including differences in dataset selection and availability of ground truth. By organizing datasets and validation approaches within the proposed taxonomy, this work contributes to a more systematic and reproducible experimental framework.

Overall, the combination of taxonomy, formalization, and evaluation-oriented analysis contributes to reducing ambiguity and improving consistency in the study of community detection methods. It also facilitates the identification of connections between approaches that are traditionally studied in isolation, such as optimization-based, probabilistic, and learning-based models.

Despite these advances, several challenges remain open and point toward promising directions for future research:

\begin{itemize}
	\item \textbf{Scalability and computational complexity:} Many high-quality methods, particularly those based on optimization or probabilistic inference, face significant scalability limitations when applied to large-scale networks.
	
	\item \textbf{Evaluation without ground truth:} In many real-world networks, true community structures are unknown or ill-defined, making validation inherently ambiguous.
	
	\item \textbf{Resolution and multi-scale structure:} Designing methods that can adapt to multiple levels of granularity without manual tuning remains an open problem.
	
	\item \textbf{Dependence on quality functions:} The strong influence of the selected modularity or objective function suggests the need for adaptive or multi-objective formulations.
	
	\item \textbf{Lack of standardized benchmarks:} The diversity of datasets and evaluation protocols hinders fair comparison across methods and calls for more unified experimental standards.
	
	\item \textbf{Integration of heterogeneous information:} Combining structural, temporal, and attribute-based information in a principled and scalable way remains challenging.
	
	\item \textbf{Robustness and adversarial settings:} Ensuring stability under noise, missing data, or adversarial manipulation is increasingly relevant.
	
	\item \textbf{Interpretability and explainability:} Bridging the gap between model complexity and interpretability is essential, especially for learning-based approaches.
\end{itemize}

In summary, this work provides a step towards a more coherent and unified understanding of the Community Detection Problem. By combining a structured taxonomy with a general mathematical formulation, it lays the groundwork for more systematic analysis, comparison, and development of future methods in this field.

{\small
\section*{Acknowledgements}
The authors acknowledge financial support by grants 
PID2020-114594GB-C21 (Mathematical optimization for machine learning in complex networks: Foundations and methodological contributions), 
PID2024-156594NB-C21 (Mathematical optimization for machine learning in complex networks: Foundations and methodological contributions), and the 
IMUS-Maria de Maeztu grant CEX2024-001517-M funded by MICIU/AEI/10.13039/501100011033.
}


\bibliographystyle{estilobibs1}
\begin{small}
\bibliography{referencias}
\end{small}


\end{document}